\documentclass[12pt]{article}

\usepackage{amsfonts,amsmath,amssymb}
\usepackage{graphicx}
\usepackage{tikz-cd}
\usepackage{enumerate}

\def\squarebox#1{\hbox to #1{\hfill\vbox to #1{\vfill}}}
\def\prod#1{\langle{#1}\rangle}

\newcommand{\qed}{\hspace*{\fill}
\vbox{\hrule\hbox{\vrule\squarebox{.667em}\vrule}\hrule}\smallskip}
\newtheorem{teorema}{Theorem}[section]
\newtheorem{lema}[teorema]{Lemma}
\newtheorem{corolario}[teorema]{Corollary}
\newtheorem{proposicao}[teorema]{Proposition}

\newenvironment{prova}{\noindent {\bf Proof:}}{\hfill $\qed $ \newline}

\newtheorem{exemplo}[teorema]{Example}

\newcommand{\R}{{\mathbb R}}

\newcommand{\Z}{{\mathbb Z}}
\newcommand{\C}{\mathbb{C}}
\newcommand{\F}{\mathbb{F}}

\renewcommand{\P}{{\mathbb P}}
\renewcommand{\t}{{\mathbb T}}

\newcommand{\ad}{{\rm ad}}
\newcommand{\Ad}{{\rm Ad}}

\newcommand{\Sl}{{\rm SL}}

\newcommand{\SU}{{\rm SU}}

\newcommand{\cl}{\mathrm{cl}}

\newcommand{\e}{{\rm e}}
\newcommand{\tr}{{\rm tr}}
\renewcommand{\max}{\mbox{{\rm max}}}
\renewcommand{\min}{\mbox{{\rm min}}}

\newcommand{\g}{\mathfrak{g}}
\renewcommand{\k}{\mathfrak{k}}
\newcommand{\s}{\mathfrak{s}}
\renewcommand{\a}{\mathfrak{a}}

\newcommand{\m}{\mathfrak{m}}

\newcommand{\n}{\mathfrak{n}}
\renewcommand{\l}{\mathfrak{l}}
\newcommand{\q}{\mathfrak{q}}
\newcommand{\p}{\mathfrak{p}}

\renewcommand{\sl}{\mathfrak{sl}}

\newcommand{\su}{\mathfrak{su}}

\newcommand{\wt}[1]{{\widetilde{#1}}}

\newcommand{\T}{\Theta}

\newcommand{\ov}[1]{{\overline{#1}}}

\renewcommand{\emptyset}{\varnothing}

\newcommand{\simto}{\stackrel{\sim}{\longrightarrow}}

\begin{document}

\title{The minimal Morse components of translations on flag manifolds are normally hyperbolic}

\author{
Mauro Patr\~{a}o\footnote{Departamento de Matem\'{a}tica,
Universidade de Bras\'{\i}lia. Campus Darcy Ribeiro.
Bras\'{\i}lia, DF, Brasil. \textit{e-mail: mpatrao@mat.unb.br }.
Supported by CNPq grant n$^{\protect\underline{\circ }}$
310790/09-3} \and
Lucas Seco\footnote{Departamento de Matem\'{a}tica,
Universidade de Bras\'{\i}lia. Campus Darcy Ribeiro.
Bras\'{\i}lia, DF, Brasil. \textit{e-mail: lseco@unb.br }. }
}

\maketitle

\begin{abstract}
Consider the iteration of an invertible matrix on the projective space: are the Morse components normally hyperbolic?
As far as we know, this was only stablished when the matrix is diagonalizable over the complex numbers. In this article we prove that this is true in the far more general context of an arbitrary element of a semisimple Lie group acting on its generalized flag manifolds: the so called translations on flag manifolds. 
This context encompasses the iteration of an invertible non-diagonazible matrix on the real or complex projective space, the classical flag manifolds of real or complex nested subspaces and also symplectic grassmanians.
Without these tools from Lie theory we do not know how to solve this problem even for the projective space.
\end{abstract}

\noindent \textit{AMS 2010 subject classification}: Primary: 37D99, 53C30 Secondary: 37B35, 22E46.

\noindent \textit{Key words:} flag manifolds, normal hyperbolicity.

\section{Introduction}

Normal hyperbolicity of an invariant manifold is the natural generalization of hyperbolicity of a fixed point, since it assures the existence of a linearization in a neighborhood of the invariant manifold \cite{pugh-shub}.
Consider the iteration of an invertible matrix on the projective space.  
%
The simplest situation is when the matrix has eigenvalues of distinct absolute values. Then the matrix is diagonalizable and the corresponding eigendirections are clearly isolated fixed points in the projective space. One can prove that they are hyperbolic fixed points and that the omega or alfa limit of any direction is one of these eigendirections.
Now if the matrix has eigenvalues with the same absolute value then the matrix is not necessarily diagonalizable and the directions in generalized eigenspaces of eigenvalues with the same absolute value give rise to a whole invariant manifold of directions. One can prove \cite{fps} that these invariant manifolds contain the alfa and omega limits of all directions and that they are the components of a minimal Morse decomposition. These Morse components are, thus, the replacement of the eigendirections. So it is natural to ask: are these Morse components normally hyperbolic?

As far as we know, this normal hyperbolicity was only stablished when the matrix is diagonalizable over the complex numbers (with the possibility of eigenvalues with the same absolute value) \cite{dkv, hermann, fps}.
In terms of the multiplicative Jordan decomposition of the matrix, previous results were not able to deal with matrices with a non-trivial unipotent component.
In this article we prove that this normal hyperbolicity is true in the far more general context of an arbitrary element of a semisimple Lie group acting on its generalized flag manifolds: the so called translations on flag manifolds. This context encompasses the iteration of an invertible non-diagonazible matrix on the real or complex projective space.  

Our approach uses techniques from Lie groups and generalized flag manifolds.
In \cite{fps, pss} we generalized \cite{dkv, hermann} and, using the Jordan decomposition, we described the Morse components in the flag manifold and their corresponding stable manifolds as orbits of certain Lie groups acting on the flag manifold. In this article we use the infinitesimal action of the Lie algebra to lift these orbits to the tangent bundle of the flag manifold and obtain natural candidates for the stable and unstable bundles of each Morse component. Then, choosing an appropriate Riemannian metric which comes from the Lie algebra, we prove the normal hyperbolicity of each Morse component.
Without these tools from Lie theory we do not know how to solve this problem even for the projective space (see Example \ref{exemplo3} and then Example \ref{exemplofinal}). 

This broader context of generalized flag manifolds encompasses other interesting cases such as the classical flag manifolds of real or complex nested subspaces and also symplectic grassmanians, which were extensively studied in the literature \cite{ammar, ayala, batterson, dkv, fps, kleinsteuber, hermann, pss, shub}. We remark that, in the wider context of flows in flag bundles, it remains an open problem to know wether the minimal Morse components are always normally hyperbolic (see \cite{pss}).


We end this introduction presenting two low dimensionsional examples of non-diagonalizable matrices where the normal hyperbolicity can be easily visualized and a third four-dimensional example where the normal hyperbolicity cannot be easily visualized.

\begin{exemplo}\label{exemplo1}
Let $X = H + N$ in $\sl(3,\R)$, where its hyperbolic and nilpotent additive Jordan components are given respectively by
\[
H = 
\begin{pmatrix}
1 & \,\, 0 & 0 \\
0 & \,\, 1 & 0 \\
0 & \,\, 0 & -2\\
\end{pmatrix}
\qquad
N = 
\begin{pmatrix}
0 & 1 & 0 \\
0 & 0 & 0 \\
0 & 0 & 0 \\
\end{pmatrix}
\]
and consider $g^t = \exp(tX)$ acting on the real projective plane $\P(\R^3)$, which is a flag manifold of the simple Lie group $G = \Sl(3,\R)$. 

    \begin{center}
        \includegraphics[scale=0.75]{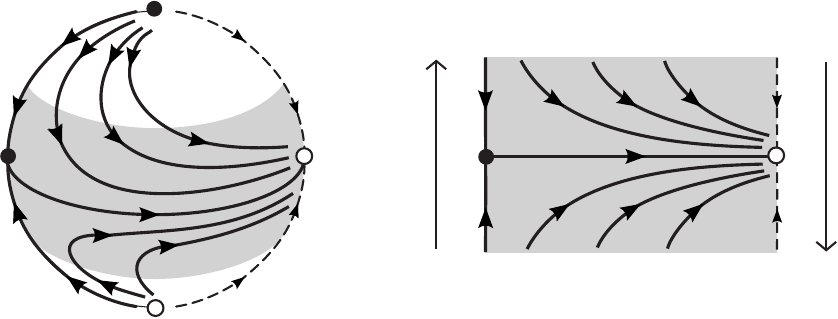}
    \end{center}

The minimal Morse components are the projectivization of the eigenspaces of $H$: the repeller is the direction $(0:0:1)$ and the attractor is the real projective line $(*:*:0)$ (a circle).
Above we sketched the phase portrait of this flow and its linearization around the attractor, which is a linear flow on the M\"obius strip over the unipotent flow on the projective line.

\end{exemplo}

\begin{exemplo}
Now, let $X = H + N$ in $\sl(2,\R) \times \sl(2,\R)$, where its hyperbolic and nilpotent additive Jordan components are given respectively by
\[
H = \left(
\begin{pmatrix}
1  & 0 \\
0  & -1 \\
\end{pmatrix}
,\,
\begin{pmatrix}
0 & 0 \\
0  & 0 \\
\end{pmatrix}
\right)
\qquad
N = \left(
\begin{pmatrix}
0 & 0 \\
0 & 0 \\
\end{pmatrix}
,\,
\begin{pmatrix}
0 & 1 \\
0 & 0 \\
\end{pmatrix}
\right)
\]
and consider $g^t = \exp(tX)$ acting on the torus, which is a flag manifold of the semisimple Lie group $G = \Sl(2,\R) \times \Sl(2,\R)$. Identifying the torus with $S^1 \times S^1$, where each $S^1$ is the projective line of $\R^2$, $g^t$ acts on the first component by the exponential of
$
\left(
\begin{smallmatrix}
1 & 0 \\
0 & -1 \\
\end{smallmatrix}
\right)
$
and acts on the second component
by the exponential of
$
\left(
\begin{smallmatrix}
0 & 1 \\
0 & 0 \\
\end{smallmatrix}
\right)
$.

    \begin{center}
        \includegraphics[scale=0.75]{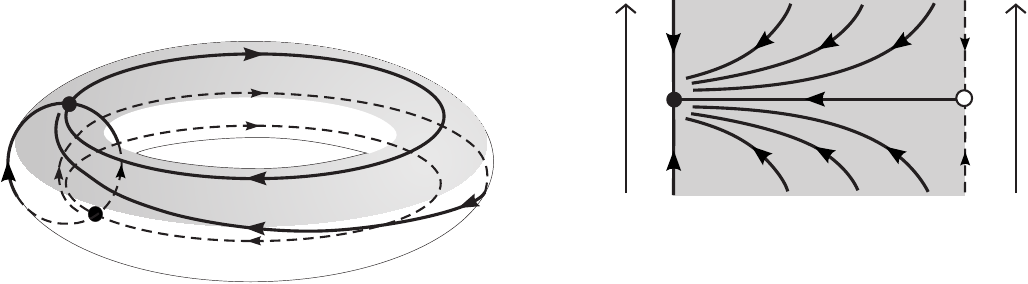}
    \end{center}

The attractor (repeller) is the cartesian product of the attractor (repeller) in the first component by $S^1$: the attractor is the $S^1$ ``above'' the torus, the repeller is the $S^1$ ``below'' the torus.
We sketched above the phase portrait of the flow on the torus and its linearization around the repeller, which is a linear flow on the cilinder over the unipotent flow on the projective line.

\end{exemplo}

\begin{exemplo}
\label{exemplo3}
Consider the same matrix $X$ of Example \ref{exemplo1} but now as element of $\sl(3,\C)$ and consider $g^t = \exp(tX)$ now acting in the complex projective space $\P(\C^3)$, which is a four dimensional real flag manifold of the simple Lie group $G = \Sl(3,\C)$. 
We have that
$$
g^t = \left( \textstyle{
\begin{array}{ccc}
e^t & te^t & 0 \\
0      & e^t  & 0 \\
0      & 0       & e^{-2t} \\
\end{array}} \right)
$$
so that 
$$
g^t (z_1: z_2: z_3) = (e^t z_1 + t e^t z_2: e^t z_2: e^{-2t} z_3)
$$
where $z_1, z_2, z_3 \in \C$ are not simultaneously zero.
If $z_3 \neq 0$ then we can assume that $z_3 = 1$ so that, dividing by $e^{-2t}$ we have
$$
g^t (z_1: z_2: 1) = (e^{3t} z_1 + t e^{3t} z_2: e^{3t} z_2: 1)
\to (0:0:1)
$$
when $t \to -\infty$ in $\R$, and tends to $(*:*:0)$ when $t \to +\infty$ in $\R$, where $\ast$ denotes arbitrary complex entries. If $z_3 = 0$ then $(z_1:z_2:0) \mapsto z_1/z_2$ is a homeomorphism with the Riemann sphere which takes $(1:0:0)$ to the infinity $\overline{\infty}$ of the complex plane, so that 
$$
g^t (z_1: z_2: 0) = (z_1 + t z_2: z_2: 0) \mapsto \frac{z_1}{z_2} + t \to \overline{\infty} = (1:0:0)
$$
when $t \to \pm\infty$ in $\R$.

The minimal Morse components are the projectivization of the eigenspaces of $H$: the attractor is the complex projective line $(*:*:0) = \P(\C^2)$ (a sphere), the repeller is the direction $(0:0:1)$.  
Since $\P(\C^3)$ is four dimensional, the stable bundle of the attractor should be a normal plane bundle over the sphere (a complex line bundle over the complex projective line) which is not easy to visualize. Thus, it is not clear if the attractor is normally hyperbolic in this case.
\end{exemplo}

The rest of the paper is organized as follows. In the second section we present the preliminaries on dynamics, homogeneous spaces of Lie groups, semi-simple Lie groups and its flag manifolds, translations on flag manifolds. In the third section we prove our main result on normal hyperbolicity.  We finish the article by revisiting Example \ref{exemplo3}.

\section{Preliminaries}\label{secpreliminar}

\subsection{Dynamics}

We recall some concepts of topological dynamics (for more
details, see \cite{alongi}). Let $\phi :\mathbb{T}\times F\to F$
be a continuous flow on a compact metric space $(F, d)$, with
discrete $\mathbb{T}={\mathbb Z}$ or continuous
$\mathbb{T}={\mathbb R}$ time. 
Denote by $\omega(x)$, $\omega^*(x)$, respectively, the forward and 
backward omega limit sets of $x$. A Morse decomposition of $\phi^t$, which is given by a finite collection of disjoint subsets $
\{{\mathcal M}_{1},\ldots ,{\mathcal M}_{n}\}$ of $F$ such that
\begin{itemize}
\item[(i)] each ${\mathcal M}_i$ is compact and
$\phi\,^t$-invariant,

\item[(ii)] for all $x \in F$ we have $\omega(x),\, \omega^*(x)
\subset \bigcup_i {\mathcal M}_i$,

\item[(iii)] if $\omega(x),\, \omega^*(x) \subset {\mathcal M}_j$
then $x \in {\mathcal M}_j$.
\end{itemize}
The minimal Morse decomposition is a Morse decomposition which is contained in every other Morse decomposition. Each set ${\mathcal M}_i$ of a minimal Morse decomposition is called a minimal Morse component. The stable/unstable set of a morse component ${\mathcal M}_i$ is the set of all points whose forward/backward omega limit set is contained in ${\mathcal M}_i$.

Now let $\phi$ be diffeomorphism on a Riemannian manifold $F$ and $D\phi$ its derivative. An invariant submanifold ${\cal M} \subset F$ is normally hyperbolic if the tangent bundle of $F$ over ${\cal M}$ has invariant vector subbundles $V^+$ and $V^-$ and positive constants $c$ and $\lambda < \mu$ such that
\begin{enumerate}[(i)]
\item $TF|_{\cal M} =  T{\cal M} \oplus V^- \oplus V^+ $
\item $|D \phi^n v| \leq c e^{-\lambda n} |v| \text{ for all } v \in V^- \text{ and } n \geq 0$
\item $|D \phi^n v| \leq c e^{\lambda n} |v| \text{ for all } v \in V^+ \text{ and } n \leq 0$
\item $|D \phi^n v| \leq c e^{\mu |n|} |v| \text{ for all } v \in T{\cal M} \text{ and } n \in \Z$
\end{enumerate}
in this case, $V^-$ is said to be the stable bundle and $V^+$ the unstable bundle of ${\cal M}$.
If $\phi^t$ is a differentiable flow on $F$, $t \in \R$, an invariant submanifold ${\cal M}$ is normally hyperbolic if its is normally hyperbolic for the time one diffeomorphism $\phi^1$.

\subsection{Homogeneous spaces of Lie groups}
\label{homogspaces}

For the theory of Lie groups and its homogeneous spaces we refer to Helgason \cite{helgason} and for the theory of principal bundles we refer to Steenrod \cite{steenrod}.  Let $G$ be a real Lie group with Lie algebra $\g$ where $g \in G$ acts on $X \in \g$ by the adjoint action
$g X = \Ad( g ) X$.  We have that  $\Ad( \exp(X) ) = e^{\ad(X)}$
where $\exp: \g \to G$ is the exponential of $G$, $\Ad$ and $\ad$ are, respectively, the adjoint representation of $G$ and $\g$.

Let a Lie group $G$ act on a manifold $F$ on the left by the differentiable map
$
G \times F \to F$,  $(g,x) \mapsto gx
$.
Fix a point $x \in F$.  The isotropy subgroup $G_{x}$ is the set of all $g \in G$ such that $gx = x$.  We say that the action is transitive or, equivalently, that $F$ is a homogeneous space of $G$, if $F$ equals the orbit $Gx$ of $x$ (and hence the orbit of every point of $F$). 
In this case, the map
$$
G \to F\qquad g \mapsto g x
$$
is a submersion onto $F$ which is a differentiable locally trivial principal fiber bundle with structure group the isotropy subgroup $G_{x}$.
Quotienting by $G_{x}$ we get the diffeomorphism
$$
G/G_{x} \simto F\qquad g G_{x} \mapsto g x
$$

If $L$ is a Lie subgroup of $G$, the orbit $L x$ is the set of all $hx$, $h \in L$. The restriction of the principal fiber bundle $G \to F$ to $L$ gives the submersion onto the orbit $L x$
$$
L \to L x \qquad l \mapsto l x
$$
which is a differentiable locally trivial principal fiber bundle with structure group $L_x = L \cap G_{x}$.  If $L x$ is an embedded submanifold of $F$ then around every point in $L x$ there exists a differentiable local section from $L x$ to $L$ that is a restriction of a local section from $F$ to $G$ of the principal fiber bundle $G \to F$.


Since the map $G \to F$ is a submersion, the derivative of the map $g \mapsto gx$ on the identity gives the infinitesimal action of $\g$, more precisely, a surjective linear map
$$
\g \to TF_x \qquad Y \mapsto Y \cdot x
$$
whose kernel is the isotropy subalgebra $\g_x$, the Lie algebra of $G_x$. 
The derivative of the map $g: F \to F$, $x \mapsto gx$, gives the action of 
$G$ on tangent vectors $g v = D g (v)$, $v \in TF$, which is related to the infinitesimal action by
$$
g( Y \cdot x ) = gY \cdot gx
$$
For a subset $ \q \subset \g$, denote by $\q \cdot x$ the set of all tangent vectors $Y \cdot x$, $Y \in \q$.
It follows that $ TF_{gx} = g( \g \cdot x ) $. In particular, for $l \in L$, the tangent space of the orbit $L x$ at $l x$ is given by $l (\l \cdot x) \subset TF_{lx}$, where $\l \subset \g$ is the Lie algebra of $L$. Thus, the tangent bundle of the orbit is given by
$$
T(L x) = L( \l \cdot x )
$$
Let $E$ be another manifold with a differentiable action of $G$, a map $f: F \to E$ is said to be $G$-equivariant if $f(g x) = g f(x)$. Such a $G$-equivariant map is automatically differentiable. 

\subsection{Semi-simple Lie theory}\label{section-lie}

For the theory of real semisimple Lie groups and their flag
manifolds we refer to Duistermat-Kolk-Varadarajan \cite{dkv},
Helgason \cite{helgason}, Knapp \cite{knapp} and Warner \cite{w}.
%
%
Let $G$ be a connected real Lie group with semi-simple Lie algebra $\g$.
Fix a Cartan decomposition $\g = \k \oplus \s$ and denote by $\langle \cdot, \cdot \rangle$ the associated Cartan inner product. 
Let $K$ be the connected subgroup with Lie algebra $\k$, it is a maximal compact subgroup of $G$.
Since $\ad(X)$ is anti-symmetric for $X \in \k$, the Cartan inner product is $K$-invariant.
Since $\ad(X)$ is symmetric for $X \in \s$, a maximal abelian subspace $\frak{a} \subset \frak{s}$ can be simultaneously diagonalized so that $\g$ splits as an orthogonal sum of
$$
\g_\alpha = \{ X \in \g:\, \ad(H)X = \alpha(H)X, \, \forall H \in \a \} 
$$
where $\alpha \in \a^*$ (the dual of $\a$). We have that $\g_0 = \m \oplus \a$, where $\m$ is the centralizer of $\a$ in $\k$. A root is a functional $\alpha \neq 0$ such that its root space $\g_\alpha \neq 0$, denote the set of roots by $\Pi$. We thus have the root space decomposition of $\g$, given by the orthogonal sum
$$
\g = \m \oplus \a \oplus \sum_{\alpha \in \Pi} \g_\alpha
$$

Fix a Weyl chamber $\frak{a}^{+}\subset \frak{a}$ and let $\Pi^{+}$ be the corresponding positive roots, $\Pi^- = - \Pi^+$ the negative roots and $\Sigma $ the set of simple roots.  Fix a subset of simple roots $\Theta \subset \Sigma $ and consider the nilpotent subalgebras
$$
\frak{n}_{\Theta }^{\pm }=\sum_{\alpha \in \Pi ^{\pm }-\langle
\Theta \rangle }\frak{g}_{\alpha }
\qquad\text{and}\qquad
\n^\pm(\T) = \sum_{\alpha \in \Pi^{\pm} \cap \langle
\Theta \rangle }\frak{g}_{\alpha }
$$
where $\langle \Theta \rangle$ is the set of roots given by linear combinations of roots in $\Theta$. Let $\n^\pm = \n^\pm_{\emptyset}$, then $\n^\pm = \n^\pm_{\T} \oplus \n^\pm(\T)$.  The minimal parabolic subalgebra is given by $\p = \m \oplus \a \oplus \n^+$ and the standard parabolic subalgebra $\frak{p}_{\Theta }$ of type $\Theta \subset \Sigma$ is given by
$$
\p_\T = \n^-(\T) \oplus \p
$$ 
so by the root space decomposition we have the orthogonal sum 
\[
\g = \frak{n}_{\Theta }^{-} \oplus \frak{p}_{\Theta }
\]
Let $p$ the dimension of $\p_\T$ and denote the grassmanian of
$p$-dimensional subspaces of $\g$ by ${\rm Gr}_p(\g)$. The flag
manifold of type $\Theta $ is the orbit
$$
\F_\T =G\frak{p}_{\Theta } \subset {\rm Gr}_p(\g)
$$
with base point $b_{\Theta }=\frak{p}_{\Theta }$ whose isotropy subalgebra is $\frak{p}_{\Theta }$ itself and isotropy subgroup is the parabolic subgroup $P_\T$.  
It follows that $\F_\T$ has dimension $\dim(\frak{n}_{\Theta }^{-})$ and that
$$
G/P_\T \simto \F_\T \qquad gP_\T \mapsto g b_\T
$$
is a $G$-equivariant diffeomorphism.  We also have that $K$ acts transitively in $\F_\T$ with isotropy subgroup $K_\T = K \cap P_\T$ so that 
$K/K_\T \simto \F_\T$, $kK_\T \mapsto k b_\T$,
is a $K$-equivariant diffeomorphism.

The Weyl group $W$ is the finite group generated by the reflections over the root hyperplanes $\alpha=0$ in $\frak{a}$, $\alpha \in \Pi$. $W$ acts on $\frak{a}$ by isometries and can be alternatively be given as $W=M^{*}/M$ where $M^{*}$ and $M$ are the normalizer and the centralizer of $\a$ in $K$,
respectively.
An element $w$ of the Weyl group $W$ can act in $\g$ by taking a representative in $M^*$.  This action centralizes $\a$, normalizes $\m$, permutes the roots $\Pi$ and thus permutes the root spaces $\g_\alpha$, where $w \g_\alpha = \g_{w \alpha}$ does not depend on the representative chosen in $M^*$.  We thus have the basepoint $w b_\T = w \p_\T$ whose isotropy subalgebra $w \p_\T$ has the orthogonal complement $w \n^-_\T$ in $\g$, that is
$$
\g = w \frak{n}_{\Theta }^{-} \oplus w \frak{p}_{\Theta }
$$

For the description of the flow $h^t = \exp (tH)$, $t \in \R$, induced by $H\in \mathrm{cl}\frak{a}^{+}$ on the flag manifold ${\mathbb F}_{\Theta }$ see (\cite{dkv}, Section 3). Its connected set of fixed points is labeled by $w \in W$, each one given by the orbit 
\[
\mathrm{fix}_{\Theta }(H,w)=G_{H} wb_{\Theta }=K_{H} wb_{\Theta },
\]
which is an embedded submanifold of $\F_\T$, where $G_{H}$ and $K_H$ denote the centralizer of $H$ respectively in $G$ and $K$ and $\mathrm{fix}_{\Theta }(H,w)$.
Consider 
the nilpotent subalgebras
\[
\frak{n}^\pm_H = \sum_{\pm \alpha(H) > 0} \g_\alpha
\]
given by the the sum of the positive/negative eigenspaces of 
$\ad(H)$ in $\g$. Since $G_H$ leaves invariant each eigenspace of $\ad(H)$ it follows that $\frak{n}^\pm_H$ is $G_H$-invariant.
Let $N_H^{\pm}$ be the corresponding connected Lie subgroups, then
$$
N_H^{\pm} \mathrm{fix}_{\Theta }(H,w)
$$
is an embedded submanifold of $\F_\T$ which is the unstable/stable manifold of $\mathrm{fix}_{\Theta }(H,w)$.

\subsection{Translations on flag manifolds}

Here we collect some previous results about the dynamics of a flow $g^t$ of translations of a real semisimple Lie group $G$ acting on its flag manifolds $\F_\T$. The flow $g^t$ is either given by the iteration of some $g \in G$, for $t \in \Z$, or by $\exp(tX)$, for $t \in \R$, where $X \in \g$, and $g^t$ acts on $\F_\T$ by left translations. Since $G$ acts on its flag manifolds by the adjoint action we will assume that $G$ is a linear Lie group, and thus $\g$ is linear Lie algebra.

The usual additive Jordan decomposition writes a matrix as a commuting sum of a semisimple and a nilpotent matrix and we can decompose the semisimple part further as the commuting sum of its imaginary and its real part, where each part commutes with the nilpotent part and the matrix is diagonalizable over the complex numbers iff its nilpotent part is zero.
This generalizes to a multiplicative Jordan decomposition of 
the flow $g^t$ in the semisimplie Lie group $G$ (see Section 2.3 of \cite{fps}), providing us with a commutative decomposition 
\[
g^t = e^t h^t u^t
\]
where there exist a Cartan decomposition of $\g$ with a corresponding maximal compact subgroup $K$ and a Weyl chamber $\a^+$ such that the elliptic component $e^t$ lies in $K$, the hyperbolic component is such that $h^t = \exp(tH)$, where $H \in \cl \a^+$, and the unipotent component is such that $u^t = \exp(tN)$, with $N \in \g$ nilpotent. Furthermore, we have that $h^t$, $e^t$ and $u^t$ lie in $G_H$, the centralizer of $H$ in $G$, and that $g^t$ is diagonalizable iff its unipotent part is $u^t = 1$.
We have that the hyperbolic component $H$ dictates the minimal Morse components (see Proposition 5.1 and Theorem 5.2 of \cite{fps}).

\begin{proposicao}
The minimal Morse components of $g^t$ on $\F_\T$ are given by ${\rm fix}_\T(H,w)$, $w \in W$, and their unstable/stable manifolds are given by $N^\pm_H {\rm fix}_\T(H,w)$. 
\end{proposicao}

\section{Normal hyperbolicity}

First we construct an appropriate Riemannian metric of $\F_\T$.  Fix the Cartan inner product $\langle \cdot, \cdot \rangle$ in $\g$ and recall that it is $K$-invariant. Let $\g_x$ be the isotropy subalgebra at $x \in \F_\T$, then
$$
\g = \g_x^\perp \oplus \g_x
$$
where $\perp$ denotes orthogonal complement with respect to the Cartan inner product. Let $k \in K$, since the isotropy subalgebra satisfies $\g_{kx} = k \g_x$, by the $K$-invariance of the Cartan inner product we have that 
$$
\g_{kx}^\perp = (k\g_x)^\perp = k( \g_x^\perp )
$$ 
Note that
\begin{equation}
\label{normalhypeq1}
\g_x^\perp \to T(\F_\T)_x \qquad X \mapsto X \cdot x
\end{equation}
is a linear isomorphism. Define the inner product in $T(\F_\T)_x$
$$
\langle X \cdot x,\, Y  \cdot x \rangle_x = \langle X, Y \rangle\qquad
\text{ where } X , Y \in \g_x^\perp
$$

\begin{proposicao}
\label{normanatural}
We have that $\langle \cdot, \cdot \rangle_x$ defines a $K$-invariant Riemannian metric of $\F_\T$ such that the map (\ref{normalhypeq1}) is an isometry. Furthermore, for $ Y \in \g $ we have
$$
| Y \cdot x |_x \leq |Y|
$$
with equality iff $Y \in \g_x^\perp$.
\end{proposicao}
\begin{prova}
Let $X \in \g_x^\perp$, then $k( X \cdot x ) = kX \cdot kx$, where $k X \in k(\g_x^\perp) = \g_{kx}^\perp$. The same holds for $kY \in \g_{kx}^\perp$, thus by the $K$-invariance of the Cartan inner product we have
$$
\langle k(X \cdot x), k(Y \cdot x) \rangle_{kx} = 
\langle kX \cdot kx, kY \cdot kx \rangle_{kx} = 
\langle kX, kY \rangle =
\langle X, Y \rangle =
\langle X \cdot x, Y  \cdot x \rangle_x
$$
To prove the smoothness of this metric, consider the local charts $\psi_s$ of $TM$ constructed from a local section $s: U \subset \F_\T \to K$ of the projection $K \to \F_\T$, $k \mapsto k b_\T$, as follows
$$
\psi_s: U \times \g_{b_\T}^\perp \to TM
\qquad
(x, Y) \mapsto s(x)( Y \cdot b_\T )
$$
Since $s(x) \in K$ and $s(x) b_\T = x$, it follows that $s(x)$ maps $\g_{b_\T}^\perp$ to $\g_x^\perp$.
By the $K$-invariance it follows that
$
\langle \psi_s(x,X), \psi_s(x,Y) \rangle_x =
\langle X, Y \rangle
$
which proves the smoothness.

For the last property, write $Y = Y_1 + Y_2$ according to $\g = \g_{x}^\perp \oplus \g_{x}$. Then $Y \cdot x = Y_1 \cdot x$, and thus
$$
|Y \cdot x|_x = |Y_1| \leq |Y|
$$
with equality iff $Y_2 = 0$ iff $Y = Y_1 \in \g_{x}^\perp$.
\end{prova}

This construction of an invariant metric is related to reductive homogenous spaces of $K$ (see \cite{helgason}) but as a model for the tangent space, instead of a subspace of the Lie algebra of $K$ as usual, here we use a subspace of the Lie algebra of $G$: this will be more appropriate to the study of $G$-action in what follows.

Now we recall the candidates for the stable and unstable vector subbundles presented in \cite{pss} which complement the tangent bundle of each Morse component ${\cal M} = \mathrm{fix}_{\T }(H,w)$.
The tangent bundle of ${\cal M} =  G_H w b_\T$ inside $T \F_\T$
is given by (see Section \ref{homogspaces})
$$
T {\cal M} = G_H ( \g_H \cdot w b_\T ) \subset T \F_\T
$$
Recall the orthogonal decomposition 
$$
\g =\g_H \oplus \n^-_H \oplus \n^+_H
$$ 
where each $\g_H$, $\n^\pm_H$ is $G_H$-invariant. Define
$$
V^\pm = G_H ( \n_H^\pm \cdot w b_\T ) \subset T \F_\T
$$

\begin{proposicao}
\label{fibradosnormais}
We have that $V^\pm$ and $T {\cal M}$ are $G_H$-invariant vector subbundles of  $T \F_\T$ over ${\cal M}$ and we have the orthogonal Whitney sum
$$
T \F_\T|_{\cal M} = T {\cal M} \oplus V^- \oplus V^+
$$
In particular $V^- \oplus V^+$ is the normal subbundle of $T {\cal M}$.
Furthermore, $v \in T {\cal M}_x$ or $v \in V^\pm_x$ can be written uniquely as $Y \cdot x$, for 
$$
Y \in \g_H \cap \g_x^\perp
\qquad
\text{or}
\qquad
Y \in \n^\pm_H \cap \g_x^\perp
$$
respectively. 
\end{proposicao}
\begin{prova}
The $G_H$-invariance is immediate from the definitions.
To prove that $V^-$ is a subbundle, first note that, since ${\cal M} = K_H w b_\T$ and $\n^-_H$ is $K_H$ invariant, it follows that $V^- = K_H ( \n^-_H \cdot w b_\T )$.  By the orthogonal decomposition
$$
\g = w \n^-_\T \oplus w \p_\T
$$
and by the root space decomposition, we have that
\begin{equation}
\label{fibradosnormaiseq1}
\n^-_H = (\n^-_H \cap w \n^-_\T) \oplus (\n^-_H \cap w \p_\T)
\end{equation}
Since we have the isotropy subalgebra $\g_{w b_\T} = w \p_\T$, it follows that
$\g_{w b_\T}^\perp = w \n^-_\T$. 
Let $x \in {\cal M}$. We have that $x = k w b_\T$, $k \in K_H$, so that 
$k \n^-_H = \n^-_H$,
$$
k w \p_\T = k \g_{w b_\T} = \g_x
\quad\text{and}\quad
k w \n^-_\T = k \g^\perp_{w b_\T} = \g^\perp_x
$$
It follows that
\begin{equation}
\label{fibradosnormaiseq2}
\n^-_H = (\n^-_H \cap \g_x^\perp) \oplus (\n^-_H \cap \g_x)
\end{equation}
so that the map
$(\n^-_H \cap \g_x^\perp) \to V^-_x$,
$Y \mapsto Y \cdot x$,
is a linear isomorphism.

Now let us prove local triviality. Since ${\cal M} = K_H w b_\T$ is an embedded submanifold of $\F_\T$, there exists a differentiable local section $\wt{s}: \wt{U} \to K_H$ of the projection $K_H \to {\cal M}$, $l \mapsto l w b_\T$, on a neighbourhood $\wt{U}$ of ${\cal M}$, such that $\wt{s}$ is the restriction of a local section ${s}: {U} \to K$ of the projection $K \to \F_\T$, $k \mapsto k w b_\T$, on neighbourhood ${U}$ of $\F_\T$. Consider the local chart $\psi_s(x,Y) = s(x)( Y \cdot w b_\T )$ of $T(\F_\T)$ as in the previous proof.  It follows that $\psi_s$ restricted to $\wt{U} \times (\n^-_H \cap 
\g_{wb_\T}^\perp
)$ is a local chart of $V^-$ given by
$$
\psi_{\wt{s}}: \, 
(x, Y) \mapsto \wt{s}(x)( Y \cdot w b_\T )
$$
since $\wt{s}(x) \in K_H$ and $\wt{s}(x) w b_\T = x$ implies that $\wt{s}(x)$ maps $\n^-_H \cap \g_{wb_\T}^\perp$ to $\n^-_H \cap \g_x^\perp$. This shows that $V^-$ is a vector subbundle.

Since equations (\ref{fibradosnormaiseq1}) and (\ref{fibradosnormaiseq2}) also hold for $\n^+_H$, the same arguments holds for $V^+$, showing that it is also a vector subbundle.  Indeed, equations (\ref{fibradosnormaiseq1}) and (\ref{fibradosnormaiseq2}) also hold for $\g^+_H$, so it follows that
$$
\g_x^\perp = (\g_H \cap \g_x^\perp) \oplus (\n^-_H \cap \g_x^\perp) \oplus (\n^+_H \cap \g_x^\perp)
$$
is an orthogonal sum and that the maps
$$
(\g_H \cap \g_x^\perp) \to T{\cal M}_x,
\qquad
(\n^\pm_H \cap \g_x^\perp) \to V^\pm_x
$$
given by $Y \mapsto Y \cdot x$ are linear isomorphisms.
In particular, since the map (\ref{normalhypeq1}) is an isometry, it follows that 
$
T \F_\T|_{\cal M} = T {\cal M} \oplus V^- \oplus V^+
$
is an orthogonal Whitney sum.
\end{prova}


These constructions of vector bundles are related to associated bundles of the principal bundle $K_H \to K_H w b_\T$ (see \cite{steenrod}) but here we constructed them as subbundles of $TM$: this will be more appropriate to the study of the $G$-action in what follows.

In order to prove normal hyperbolicity we need the following lemma.

\begin{lema}\label{lemadecaimentoexp}
Let $H \neq 0$. We have that
\[
|h^t Y| \leq {\rm e}^{-\mu t} |Y|,\quad\mbox{for}\quad Y \in
\n^-_H,\quad t \geq 0
\]
and
\[
|h^t Y| \leq {\rm e}^{\mu t} |Y|,\quad\mbox{for}\quad Y \in
\n_H,\quad t \leq 0
\]
where
\[
\mu = \min\{\alpha(H) :  \alpha(H)>0, \, \alpha \in \Pi\}
\]
\end{lema}
\begin{prova}
For $Y \in \n^\pm_H$, we have that $h^t Y = \e^{t\ad(H)} Y$, where
$\e^{t\ad(H)}$ is $\prod{\cdot,\cdot}$-symmetric with eigenvalues in $\n^\pm_H$ given by
\[
\{ \e^{\pm\alpha(H) t}:\, \quad \alpha(H) > 0,\, \alpha \in \Pi\}
\]
since $\ad(H)$ is $\prod{\cdot,\cdot}$-symmetric with eigenvalues in $\n^\pm_H$ given by
\[
\{ \pm\alpha(H):\, \quad \alpha(H) > 0,\, \alpha \in \Pi\}
\]
Writing $Y$ as the orthogonal sum of eigenvectors $Y = \sum_\alpha
Y_\alpha$, we have that
\[
|h^t Y| = |\sum_\alpha e^{\pm\alpha(H) t} Y_\alpha | \leq
\sum_\alpha e^{\pm\alpha(H) t} | Y_\alpha |
\]
For $t > 0$ and $Y \in \n^-_H$, we have that
\[
|h^t Y| \leq e^{-\mu t} \sum_\alpha | Y_\alpha | = e^{-\mu t}|Y|
\]
since $e^{-\alpha(H) t} < e^{-\mu t}$, for all $\alpha \in
\Pi$ with $\alpha(H) > 0$. For $t < 0$ and $Y \in \n_H$, we have that
\[
|h^t Y| \leq e^{\mu t} \sum_\alpha | Y_\alpha | = e^{\mu t}|Y|
\]
since $e^{\alpha(H) t} < e^{\mu t}$, for all $\alpha \in
\Pi$ with $\alpha(H) > 0$.
\end{prova}

We now prove normal hyperbolicity. 

\begin{teorema}\label{teo}
Each Morse component $\mathrm{fix}_{\T }(H,w)$ is normally hyperbolic.
\end{teorema}
\begin{prova}
By the Jordan decomposition of $g^t$ we have the following commutative decompostition
\[
g^t = e^t h^t u^t
\]
where $h^t = \exp(tH)$, with $H \in \g$ hyperbolic, $u^t = \exp(tN)$, with $N \in \g$ nilpotent, and $e^t, u^t \in G_H$, the centralizer of $H$ in $G$. Furthermore, we can assume that $H \in \cl \a^+$ and that $e^t \in K_H$, the centralizer of $H$ in $K$.

By Proposition \ref{fibradosnormais} a tangent vector $v \in V^\pm_x$ can be written as $v = Y \cdot x$, for $Y \in \n^\pm_H \cap \g_x^\perp$.  
By Proposition \ref{normanatural} we have that $|v| = |Y|$ and also that
$$
|g^t v| = |g^tY \cdot g^t x| \leq |g^tY|
$$
Since $g^t \in G_H$ implies $g^t Y \in \n_H^\pm$, it is enough to show that the inequalities hold for $g^t$ restricted to $\n_H^\pm$.
This follows from standard linear algebra and we will sketch the argument here for the readers' convenience.  By Lemma \ref{lemadecaimentoexp}, there exists $\mu >0$ such that $|h^t X| \leq \e^{-\mu t} |X|$, for $t \geq 0$ and $X \in \n_H^-$.  Since we can assume that $e^t \in K_H$, it follows that
\[
|g^tY| = |h^tu^tY| \leq \e^{-\mu t} |u^t Y|
\]
Since $u^t = \exp(t N)$, for some nilpotent $N \in \g$, we have that $u^t Y = \e^{t\ad(N)} Y$. By the triangle inequality, we have that
\[
|u^t Y| = |\e^{t\ad(N)} Y| \leq \sum_{k \geq 0} \frac{t^k}{k!} \| \ad(N)^k
\| |Y| = p(t) |Y|
\]
where $\| \cdot \|$ is the operator norm induced by the norm $|\cdot |$ in $\g$ and $p(t)$ is a polynomial, since $\ad(N)$ is nilpotent.  Thus, for $v \in V^-$, we have that
\[
|g^t v| \leq \e^{-\mu t} p(t) |v|, \quad t > 0
\]
The case for $V^+$ is completely analogous and we get
\[
|g^t v| \leq \e^{\mu t} p(t) |v|, \quad t < 0
\]

For $T {\cal M}$, note that for $x \in {\cal M} = \mathrm{fix}_{\T }(H,w)$ we have $g^t x = e^t u^t x$, thus $g^t$ acts as $e^t u^t$ in $T {\cal M}$.  By Proposition \ref{fibradosnormais} a tangent vector $v \in T{\cal M}_x$ can be written as $v = Y \cdot x$, for $Y \in \g_H \cap \g_x^\perp$. By Proposition \ref{normanatural} we have that $|v| = |Y|$ and also that 
$$
|g^t v| = |e^t u^t Y \cdot e^t u^t x| \leq |e^t u^t Y| 
= |u^t Y| \leq p(t) |Y| = p(t) |v|
$$
where we used that $e^t \in K_H$ and the same inequality for $|u^t Y|$ as above.

Since $\e^{-\frac{\mu}{2} t} p(t) \to 0$ as $t \to +\infty$, it is bounded by $c_1$ for $t \geq 0$, so that
$$
\e^{-\mu t} p(t) = \e^{-\frac{\mu}{2} t} \left( \e^{-\frac{\mu}{2} t} p(t) \right)
\leq c_1  \e^{-\frac{\mu}{2} t}, \quad t \geq 0
$$
By the same argument
$$
\e^{\mu t} p(t) = \e^{\frac{\mu}{2} t} 
\left( \e^{\frac{\mu}{2} t} p(t) \right)
\leq c_2  \e^{\frac{\mu}{2} t}, \quad t \leq 0
$$
At last, since $\e^{ -\mu |t|} p(t) \to 0$ as $t \to \pm\infty$, it is bounded by $c_3$ for $t \in \t$, so that
$$
p(t) = \e^{\mu |t|} \left( \e^{ -\mu |t|} p(t) \right)
\leq c_3 \e^{\mu |t|}, \quad t \in \t
$$
Conditions (ii), (iii) and (iv) of normal hyperbolicity then follows by choosing $\lambda = \frac{\mu}{2}$
and $c = \max\{ c_1, c_2, c_3 \}$.  
\end{prova}

It follows that $V^\pm$ are the unstable/stable bundle of $g^t$.
By the main result of \cite{pugh-shub}, we obtain a linearization of this flow in a neighborhood of each minimal Morse component ${\cal M} = {\rm
fix}_\T(H,w)$.

\begin{corolario}\label{corlinearizgt}
Let $V = V^- \oplus V^+$.
There exists a differentiable map $ V \to
{\mathbb F}_{\T }$ which takes the null section to 
${\rm fix}_\T(H,w)$ and such that:

\begin{enumerate}[$(i)$]
\item Its restriction to some neighborhood of the null section $V_0$ inside $V$ is a $g^t$-equivariant diffeomorphism onto some neighborhood of ${\rm fix}_\T(H,w)$ inside $\F_\T$.

\item Its restrictions to $V^\pm$ are $g^t$-equivariant diffeomorphisms, respectively, onto the unstable/stable manifolds 
$N^\pm_H \mathrm{fix}_{\T}(H,w)$.
\end{enumerate}
\end{corolario}
\begin{prova}
It is enough to note that the action of $g^t$ on $V$ is given by the restriction of the differential of the action of $g^t$ on $\F_\T$ and also that the equivariance property is equivalent to the conjugation property of \cite{pugh-shub}.
\end{prova}

We remark that, in the wider context of flows in flag bundles, it remains an open problem to know wether the minimal Morse components are always normally hyperbolic.  In \cite{pss} our main tool requires a $G_H$-equivariant linearization of the flow generated by $H$ around a connected component of its fixed point set on the flag manifold which, unfortunately, we were only able to construct in some situations.  It would be nice if the linearization we get in the previous result could be made $G_H$-equivariant so that it could be used to provide the linearization of the flows on flag bundles.

We end this article by revisiting Example \ref{exemplo3} of the introduction.

\begin{exemplo}
\label{exemplofinal}
Given a non-null vector $v \in \C^n$, denote its corresponding direction by $[v] = \C v$. Then the complex projective space $\P(\C^n)$ is the set of such directions and the action of an invertible $n \times n$ matrix $g$ on $\P(\C^n)$ is given by $g [v] = [g v ]$.  If $v = (z_1,\ldots,z_n)$ we denote its direction by $[v] = (z_1:\ldots:z_n)$.  The canonical complex line bundle $\gamma(\C^n)$ over $\P(\C^n)$ is the vector bundle given by
$
\gamma(\C^n) = \{ (x, v):\,  v \in x \} \subset \P(\C^n) \times \C^n
$.

On $\P(\C^3)$, consider the flow $g^t = \exp(tX)$ of translations given in Example \ref{exemplo3}, where the matrix $X$ is  an element of $\sl(3,\C)$ and has hyperbolic component
$$
H = \left(
\begin{smallmatrix}
1 & \,\,   &   \\
  & \,\, 1 &   \\
  & \,\,   & -2\\
\end{smallmatrix} 
\right)
$$
We have that $\P(\C^3)$ is a flag manifold of the simple real Lie group $G = \Sl(3,\C)$ of complex matrices of determinant one, whose Lie algebra is the real Lie algebra $\g = \sl(3,\C)$ of complex traceless matrices.
In order to see this, let $\k = \su(3)$ be the subalgebra of anti-hermitian matrices and let $\s$ be the subspace of hermitian matrices in $\g$. This gives a Cartan decomposition $\g = \k \oplus \s$ whose corresponding Cartan inner product is proportional to the usual hermitian inner product on complex matrices, more precisely we have 
$
\langle X, Y \rangle = 6 \, \tr( X \ov{Y}^\dagger )
$,
where $\ov{Y}^\dagger$ is the conjugate transpose of $Y$.
The subgroup $K = \SU(3)$ of unitary matrices $k$ of determinant one, that is, $k \ov{k}^\dagger = I$ and $\det(k)=1$, has Lie algebra $\k$ and is a maximal compact subgroup of $G$.  The subspace $\a \subset \s$ of real diagonal matrices is a maximal abelian subalgebra of $\s$, the corresponding root spaces are given by $\C E_{ij}$, $i \neq j$, where $E_{ij}$ is the elementary matrix with entry $1$ in row $i$ column $j$ and entries $0$ elsewhere. 
An element $H = {\rm diag}( \lambda_1, \lambda_2, \lambda_3) \in \a$ is such that $\ad(H)$ has eigenvalue $\alpha_{ij}(H) = \lambda_i - \lambda_j$ in $\C E_{ij}$. 
Fix the Weyl chamber given by the subset $\a^+ \subset \a$ of diagonal matrices with decreasing diagonal entries.
The positive root spaces are $\C E_{12}$, $\C E_{13}$, $\C E_{23}$, so that
$$
\n^+ = \left\{ 
\left(
\begin{smallmatrix}
 & * & *\\
 &  & *\\
\phantom{*} &  & \\
\end{smallmatrix}
\right)
\right\}
\qquad
\text{and}
\qquad
\n^- = \left\{ 
\left(
\begin{smallmatrix}
  &   & \phantom{*} \\
* &   &  \\
* & * &  \\
\end{smallmatrix}
\right)
\right\}
$$
where $\ast$ denotes arbitrary complex entries. Since the centralizer of $\a$ in $\k$ is given by $\m = \sqrt{-1}\,\a$, we have that $\m \oplus \a = \mathfrak{d}$ the subspace of complex traceless diagonal matrices.  
It follows that the minimal parabolic subalgebra $\p = \mathfrak{d} \oplus \n^+$ is the set of upper triangular matrices in $\g$.
Fix the simple root $\T = \{ \alpha_{23} \}$ and thus the corresponding negative root space $\C E_{32}$. It follows that $\p_\T = \C E_{32} \oplus \p$ so that its corresponding parabolic subgroup 
$P_\T$ and orthogonal complement $\n^-_\T$ are given by
$$
P_\T =
\left\{
\left(
\begin{smallmatrix}
* & * & *\\
 & * & *\\
\phantom{*} & \circledast & * \\
\end{smallmatrix}
\right)
\in G
\right\}
\qquad
\text{and}
\qquad
\n^-_\T =
\left\{
\left(
\begin{smallmatrix}
\phantom{*} &   &  \\
* &   &  \\
* & \phantom{*} & \phantom{*} \\
\end{smallmatrix}
\right)
\right\}
$$
Note that $P_\T$ is precisely the isotropy subgroup of the direction $(1:0:0)$ in $\P(\C^3)$. It follows that $\P(\C^3)$ is the flag manifold of $G$ of type $\T$, more precisely, the map
$$
\F_\T \to \P(\C^3) \qquad g b_\T \mapsto g (1:0:0), \quad g \in G
$$  
is a $G$-equivariant diffeomorphism.

To describe the dynamics first note that, since the adjoint action of $G$ in $\g$ is by conjugation, by the block form of $H$ we have that $g \in G$ centralizes $H$ iff it leaves invariant its eigenspaces. Thus we have
$$
G_H =
\left\{
\left(
\begin{smallmatrix}
* & * &  \\
* & * &  \\
  &   & * \\
\end{smallmatrix}
\right)
\in G
\right\}
\quad
\text{and}
\quad
K_H =
\left\{
\left(
\begin{smallmatrix}
a & -\ov{b} &  \\
b & \phantom{-}\ov{a}  \\
  &   & 1 \\
\end{smallmatrix}
\right)
:
\begin{array}{l}
a, b \in \C \\
|a|^2 + |b|^2 = 1
\end{array}
\right\} \simeq \SU(2)
\qquad
$$
Then note that the positive eigenspaces of $\ad(H)$ are $\C E_{13}$ and $\C E_{23}$ so that
$$
\n^+_H = \left\{ 
\left(
\begin{smallmatrix}
 &   & *\\
 &   & *\\
\phantom{*} & \phantom{*} & \\
\end{smallmatrix}
\right)
\right\}
\qquad
\text{and}
\qquad
\n^-_H = \left\{ 
\left(
\begin{smallmatrix}
  &   & \phantom{*} \\
  &   &  \\
* & * &  \\
\end{smallmatrix}
\right)
\right\}
$$
The attractor is given by $w=1$ so that it is
\begin{eqnarray*}
{\cal M} & = & G_H (1:0:0) \\
& = & (*:*:0) = \P( \C^2 ) 
\end{eqnarray*}
Put $b_\T = (1:0:0)$, the stable bundle of the attractor is the normal bundle $V^- = G_H( \n^-_H \cdot b_\T )$ so that we have the orthogonal Whitney sum
$$
T \P(\C^3)|_{\P(\C^2)} = T \P(\C^2) \oplus V^-
$$
We claim that $V^-$ is the canonical complex line bundle $\gamma = \gamma(\C^2)$ over the complex projective line $\P( \C^2 )$, in particular $V^-$ is a non-trivial vector bundle.

Indeed, we have that
$$
\gamma = \left\{  ( (a:b), \, c (a,b) ):\,
\begin{array}{r}
a,b,c \in \C \\
|a|^2 + |b|^2 = 1
\end{array} \right\}
\subset \P( \C^2 ) \times \C^2
$$
By the proof of Proposition \ref{fibradosnormais}, we have that
$V^- = K_H( \n^-_H \cap \n^-_\T \cdot b_\T )$,
where
$$
\n^-_H \cap \n^-_\T =
\left\{ 
\left(
\begin{smallmatrix}
  &   & \phantom{*} \\
  &   &  \\
* & * &  \\
\end{smallmatrix}
\right)
\right\}
\cap
\left\{
\left(
\begin{smallmatrix}
\phantom{*} &   &  \\
* &   &  \\
* & \phantom{*} & \phantom{*} \\
\end{smallmatrix}
\right)
\right\}
=
\left\{ 
\left(
\begin{smallmatrix}
  &   & \phantom{*} \\
  &   &  \\
* & \phantom{*} &  \\
\end{smallmatrix}
\right)
\right\}
=
\C Y^-
$$
where we let $Y^- = E_{31}$.  Define the $K_H$-equivariant map
$$
V^- \to \gamma \qquad 
k( y Y^- \cdot b_\T ) \mapsto 
( k(1:0), \, \ov{y} k(1,0) )
$$
where $k \in K_H$ and $y \in \C$. By fixing $k$, we see clearly that it is an isomorphism between fibers thus, in order to show it is a bundle homeomorphism, we must only show it is well defined.
If $k( y Y^- \cdot b_\T )  = k'( y' Y^- \cdot b_\T )$ then $k b_\T = k' b_\T$ so that $k' = kl$ where
$$
l \in K_H \cap P_\T = 
\left\{
\left(
\begin{smallmatrix}
c &  &   \\
  & \ov{c} &   \\
  &   &  1 \\
\end{smallmatrix}
\right):\,
\begin{array}{r}
c \in \C \\
|c| = 1
\end{array}
\right\}
\simeq S^1
$$
%
%
Thus $y Y^- = l( y' Y^- ) = y' l Y^-$, where the action of $l$ on $Y^-$ is by conjugation
$$
l Y^- l^{-1} = 
\left(
\begin{smallmatrix}
c &  &   \\
  & \ov{c} &   \\
  &   &  1 \\
\end{smallmatrix}
\right)
\left(
\begin{smallmatrix}
\phantom{c} & \phantom{c} & \phantom{c} \\
& & \\
1 & & \\
\end{smallmatrix}
\right)
\left(
\begin{smallmatrix}
\ov{c} &  &   \\
  & c &   \\
  &   &  1 \\
\end{smallmatrix}
\right)
=
\left(
\begin{smallmatrix}
\phantom{c} & \phantom{c} & \phantom{c} \\
& & \\
\ov{c} & & \\
\end{smallmatrix}
\right)
=
\ov{c} Y^-
$$
It follows that $y Y^- = y' \ov{c} Y^-$, so that $y' = y c$, since $\ov{c}c = 1$. We also have that
$$
k' = kl = 
\left(
\begin{smallmatrix}
a & -\ov{b} & \\
b & \phantom{-}\ov{a} \\
  & & 1 \\
\end{smallmatrix}
\right)
\left(
\begin{smallmatrix}
c &  &   \\
  & \ov{c} &   \\
  &   &  1 \\
\end{smallmatrix}
\right)
=
\left(
\begin{smallmatrix}
ca & -\ov{cb} & \\
cb & \phantom{-}\ov{ca} \\
  & & 1 \\
\end{smallmatrix}
\right)
$$
Thus
$$
k( y Y^- \cdot b_\T ) \mapsto ( k(1:0), \, \ov{y} k(1,0) ) 
= ( (a:b), \, \ov{y} (a,b) )
$$
and
$$
k'( y' Y^- \cdot b_\T ) \mapsto 
( (ca:cb), \, \ov{yc} (ca,cb) )
=( (a:b), \, \ov{y} (a,b) )
$$
which shows the map is well defined and thus gives a vector bundle isomorphism from $V^-$ to $\gamma$.
Since $V^-$ is a vector bundle over a sphere, this fact can also be seen by computing its clutching function over the equator: it turns out to be homotopic to $z \mapsto z$, $|z|=1$, which is the clutching function of the canonical complex line bundle over $\P(\C^2)$ (see \cite{hatcher-kt}, Proposition 1.11). The real counterpart of this is Example \ref{exemplo1}, where the stable bundle of the attractor is the infinite M\"oebius strip, which is the canonical real line bundle over the real projective line.

Note that the $K$-invariant metric restricted to $V^-$ is proportional to the natural metric of $\gamma \subset \P(\C^2) \times \C$, given by the hermitian inner product in the second factor. Nevertheless, the above vector bundle isomorphism $V^- \simeq \gamma$ is only $K_H$-equivariant. Since $g^t \in G_H$ it follows that, even for this example in $\P(\C^3)$, $\gamma$ cannot be used directly to study the normal dynamics of $g^t$ over the attractor: the invariant metric and the subbundle $V^-$ constructed with tools from Lie theory are much more appropriate for this.

\end{exemplo}

\end{document}